\documentstyle[12pt,amssymb]{amsart}

\theoremstyle{plain}
\newtheorem{thm}{Theorem}
\newtheorem{theorem}[thm]{Theorem}
\newtheorem{corollary}[thm]{Corollary}
\newtheorem{lemma}[thm]{Lemma}
\newtheorem{remark}[thm]{Remark}
\newtheorem{example}[thm]{Example}

\newenvironment{theorem*}[1]{\smallskip\noindent{\bf #1.}\rm}{\medskip}
\newenvironment{proofof}[1]{\smallskip\noindent{\it #1}\rm}
                {\hspace*{\fill} $\Box$\medskip}
\newenvironment{proof}{\smallskip\noindent{\it Proof.}\rm}
                        {\hspace*{\fill} $\Box$\medskip}

\newcommand\esssup{\operatornamewithlimits{ess\,sup}}
\newcommand\spr{\operatorname{r\,}}
\newcommand\supp{\operatorname{supp\,}}
\newcommand\cs{\curlyeqsucc}
\newcommand\cp{\curlyeqprec}

\newcommand\al{\alpha}

\newcommand\de{\delta}
\newcommand\la{\lambda}
\newcommand\om{\omega}
\newcommand\eps{\varepsilon}

\newcommand\cA{{\mathcal A}}
\newcommand\cB{{\mathcal B}}
\newcommand\cH{{\mathcal H}}
\newcommand\cM{{\mathcal M}}
\newcommand\cX{{\mathcal X}}

\newcommand\bN{{\mathbb N}}

\newcommand\op{operator}

\title{On similarity of perturbed multiplication operators}
\author{R.~O.~Hryniv and Ya.~V.~Mykytyuk}

\address{Institute for Applied Problems of Mechanics and Mathematics,
3b~Naukova str., 79601 Lviv, Ukraine}
\email{hryniv@@mebm.lviv.ua}

\address{Lviv National University, 1 Universytetska str., 79602 Lviv, Ukraine}
\email{yamykytyuk@@yahoo.com}


\date{\today}

\begin{document}

\begin{abstract}
Let $S$ be the multiplication operator by an independent variable $x$ in
$L_2(0,1)$ and $V$ be an integral operator of Volterra type. We find
conditions for $T:=S+V$ to be similar to $S$ and discuss some
generalisations of the results obtained to an abstract setting.
\end{abstract}

\maketitle

\section{Introduction}

In the Hilbert space $\cH=L_2(0,1)$, consider the operator $S$
of multiplication by an independent variable,
\(
        (Sf)(x) = x f(x)
\),
and its perturbation $T:= S + V$, where
\begin{equation}\label{eq:V}
        (Vf)(x) := \int_0^x v (x,t) f(t)\, dt.
\end{equation}
We assume throughout the paper that the kernel $v$ is Lebesgue
measurable on $[0,1]\times[0,1]$, $v(x,t)=0$ if $x < t$,
and that the induced integral \op\ $V$ is bounded in~$\cH$.

Operators of similar type appear, e.g., in the so-called Friedrichs
model~\cite{Fri,Fa} or in polymerisation chemistry, where $T$ describes the
evolution of a polymer system near dynamical equilibrium~\cite{Ko}. In both
examples the asymptotic behaviour of the group~$e^{itT}$ (in particular,
uniform boundedness, or Lyapunov stability, of~$e^{itT}$) is of much
importance, which poses a problem of similarity of~$T$ to a selfadjoint
\op. For the case when $v(x,t) = \phi(x)\psi(t)$ for $0 \le t\le x \le1$ 
and $\phi,\psi \in \cH$ this problem was studied in detail in the 
paper~\cite{NT}, where the following result was established.

\begin{theorem*}{Theorem A} \em
Let $\phi\psi \equiv0$ and suppose that there exist moduli of continuity
$\om_1,\om_2$ such that $\phi\in Lip\,(\om_1)$, $\psi\in Lip\,(\om_2)$, and
\begin{equation}\label{eq:modulus}
        \int_0 \frac{\om_1(\tau)\om_2(\tau)}{\tau} <\infty.
\end{equation}
Then the operator $T$ is similar to a selfadjoint one.
\end{theorem*}

\noindent
Moreover, a sharp analysis of behaviour of the \op~$T$ resolvent near the 
real axis shows the necessity of condition~\eqref{eq:modulus} in the sense
that if it does not hold, then the \op~$T$ need not be similar to a
selfadjoint one; in~\cite{NT} the corresponding examples are constructed for 
$\om_1(\tau) = |\ln \tau|^{-\delta_1}$ and
$\om_2(\tau) = |\ln \tau|^{-\delta_2}$ with $\delta_1 + \delta_2 < 1$. 

The main aim of this paper is to find conditions under which the perturbed
\op~$T$ is similar to the unperturbed \op~$S$. We consider the 
perturbations~$V$ of the general form~\eqref{eq:V} and follow a line of 
attack due to Friedrichs~\cite[Ch.~II.6]{Fri}. Namely, we find
sufficient conditions for existence of a bounded \op~$K$ with spectral
radius $\spr(K)$ smaller than~$1$ such that
\begin{equation}\label{eq:simil}
        T (I + K) = (I + K) S.
\end{equation}
Our main results are as follows.

\begin{theorem}\label{thm:main}
 Suppose that the kernel
 \[
   w (x,t) := \frac{|v(x,t)|}{x-t}, \qquad x,t\in[0,1],
 \]
 generates an integral \op~$W$ that is bounded in~$\cH$ and has spectral
 radius $\spr(W)$ less than $1/2$. Then the \op s $T$ and $S$ are similar.
\end{theorem}

\begin{corollary}\label{cor:toepl}
 Suppose that there exists a function $q\in L_1(0,1)$ such that
 $w(x,t) \le q(x-t)$ for all $x,t,\ 0\le t<x\le 1$. Then the \op s $T$ and
 $S$ are similar.
\end{corollary}

Note that under the assumptions of Theorem~A we have
\begin{align*}
  |v(x,t)| &\le |\phi(x)\psi(t)| + |\phi(t)\psi(x)| =
        \bigl(|\phi(x)| - |\phi(t)|\bigr)\bigl(|\psi(t)| - |\psi(x)|\bigr)\\
        &\le \om_1(x-t)\om_2(x-t),
\end{align*}
and therefore Corollary~\ref{cor:toepl} applies with $q(\tau) =
\om_1(\tau)\om_2(\tau)/\tau$ and proves the claim of Theorem~A.
Corollary~\ref{cor:toepl} also admits kernels $v$ for which the norm
\[
	\|v\|_\al := \esssup_{0\le t<x\le1} (x-t)^{1-\al}|v(x,t)|
\]
is finite for some $\al>0$; see~\cite{Fre,Ma} and references therein for 
related details on similarity of Volterra operators to fractional integration 
operators. We remark that the results of~\cite{NT} imply that the condition
$q \in L_1(0,1)$ cannot be weakened.

In fact, under our approach $V$ need not be an integral \op\ of the
form~\eqref{eq:V}. We allow $V$ from the algebra~$\cA$ of operators leaving
invariant functions with support in $[a,1]$, for any $a\in [0,1)$, that are
majorised in a certain sense. The corresponding concepts are based on the
theory of nonnegative operators in a Banach space with a positive
cone~\cite{KR} and are developed in Section~\ref{sec:nonneg}. Abstract results
on similarity of the \op s $T$ and $S$ are established in Section~\ref{sec:simil},
and then used to prove Theorem~\ref{thm:main} and Corollary~\ref{cor:toepl} in
Section~\ref{sec:appl}. Finally, in the last section we comment on some
straightforward generalisations of the main results.

Throughout the paper we shall denote by $\spr(T)$ the spectral radius of a 
bounded \op~$T$; recall that 
\(
  \spr(T) = \lim_{n\to \infty} \|T^n\|^{1/n}.   
\)


\section{Nonnegative \op s and some auxiliary results}\label{sec:nonneg}

We start by recalling some concepts of linear spaces with a positive
cone (see, e.g., \cite{KR}). Denote by
\[
        \cH_+ = \{ f \in \cH \mid f(x) \ge 0 \ \text{a.e.~in } [0,1]\}
\]
and
\[
        \cB_+(\cH) = \{ A \in \cB(\cH) \mid A \cH_+ \subset \cH_+ \}
\]
the cones of nonnegative elements in~$\cH$ and nonnegative \op s in~$\cH$,
respectively. As usual, for any $f,g \in \cH$ and $A,B \in \cB(\cH)$
we write $f \ge g$ and $A\ge B$ if $f-g\in\cH_+$ and $A-B\in \cB_+(\cH)$,
respectively. Recall that $\cH_+$ is a generating cone and hence for any
$f\in\cH$ the absolute value $|f|$ exists as an element of $\cH_+$; in the
present context $|f|$ is the function defined by $|f|(x) = |f(x)|$.
The cone $\cB_+(\cH)$, on the contrary, is not generating and hence the
absolute value $|A|$ can not be defined for all $A \in \cB(\cH)$.
We shall point out the class of operators, for which the absolute value is
well defined.

An \op\ $B$ is said to {\em majorise} $A$ (written $B \cs A$ or $A \cp B$)
if $|Af| \le B|f|$ for all $f\in \cH$. Evidently, $A \cp B$ implies
that $B$ is nonnegative, $B \ge \la A$ for all $\la$ with $|\la|\le1$, and
that
\[
    |(Af,g)| \le (|Af|,|g|) \le (B|f|,|g|) \le \|B\| \|f\| \|g\|
\]
for all $f,g \in \cH$, i.e., $\|A\| \le \|B\|$. Moreover, if $B$ majorises
$A$, then the inequality
\[
        |A^nf| \le B|A^{n-1}f| \le \dots \le B^{n-1}|Af| \le B^n |f|
        \qquad \forall  f\in\cH
\]
shows that $A^n \cp B^n$ for any $n\in \bN$ and hence $\spr(A) \le \spr(B)$.

Put
\[
   \cM(A) := \{ B\in\cB_+(\cH) \mid B \cs A \}
\]
and
\[
   \cB_{\cM}(\cH) := \{ A \in \cB(\cH) \mid \cM(A) \ne \emptyset \}.
\]
Observe that $\cB_{\cM}(\cH)$ is a multiplicative cone in $\cB(\cH)$, i.e.,
$AB\in \cB_\cM(\cH)$ whenever $A,B \in \cB_\cM(\cH)$.

\begin{example}\label{eq:1}
 Let $A$ and $B$ be integral \op s in~$\cH$ with continuous kernels $a$ and
 $b$ respectively. Then $B \in \cM(A)$ if and only if $|a(x,t)| \le b(x,t)$
 everywhere in $[0,1]\times[0,1]$ and $A \in \cB_\cM(\cH)$ if and only if the
 kernel $|a|$ induces a bounded integral \op\ in~$\cH$.
\end{example}

To some extent, the above example is generic as any bounded \op\ in~$\cH$
is a strong limit of integral \op s with continuous kernels.  In
fact, let $r$ be a smooth nonnegative function such that
$\supp r \subset [0,1]$ and $\int r = 1$.
Denote by $R_\eps$ the integral operator
\[
        (R_\eps f)(x) = \int_0^1 r_\eps(x-t) f(t)\, dt,
        \qquad r_\eps(x) = \eps^{-1}r(x/\eps),
\]
 and for any $A \in \cB(\cH)$ and $\eps>0$ put
\[
        A_\eps:= R_\eps A R_\eps.
\]
 Note that $R_\eps \to I$ as $\eps \to 0$ strongly in
 $\cH$~\cite[Ch.~III.2]{St}, whence
 \[
       \text{s\,-}\lim_{\eps\to0} A_\eps = A.
 \]
 Moreover, since $R_\eps$ is a Hilbert-Schmidt \op, $A_\eps$ is an
 integral \op\ with some kernel $a_\eps$ (see details in \cite{HS}). Denoting
 by $\de_\xi(u):= \de ( u - \xi)$ the delta function centered at the
 point~$\xi$, we see that
 \[
   a_\eps(x,t) = (A_\eps \de_t,\de_x) = (AR_\eps\de_t,R^*_\eps\de_x)
        = ( A r_{\eps,t}, \hat r_{\eps,x} ),
 \]
 where $r_{\eps,t}(u) := r_\eps(u-t)$ and $\hat r_{\eps,x}(u):= r_\eps(x-t)$
 are continuous functions in $t$ and $x$ respectively with values in~$\cH$,
 whence the kernel~$a_\eps$ is continuous.  

 This observation is heavily used to justify the following statement.

\begin{lemma}\label{lem:abs}
 For any $A \in \cB_{\cM}(\cH)$, the set $\cM(A)$ contains the minimal
 element $|A|$ called the {\em absolute value} of $A$. In other words, $|A|
 \in \cM(A)$ and $|A| \le B$ for any $B \in \cM(A)$.
\end{lemma}

\begin{proof}
 Suppose that $A \in \cB_{\cM}(\cH)$ and $B \in \cM(A)$ and put
 $A_\eps = R_\eps A R_\eps$, $B_\eps:= R_\eps B R_\eps$ with the
 above constructed $R_\eps$. Then
 \begin{align*}
   |(A_\eps f,g)| = |(AR_\eps f, R^*_\eps g)|
        &\le (B|R_\eps f|, |R^*_\eps g|) \\
        &\le (BR_\eps |f|, R^*_\eps |g|) = (B_\eps |f|,|g|)
 \end{align*}
 for any $f,g \in \cH$, which implies the inequality
 \begin{equation}\label{eq:ab}
        |a_\eps(x,t)| \le b_\eps(x,t)
 \end{equation}
 for the corresponding kernels.

 Denote by $A^+_\eps$ the integral \op\ induced by the kernel
 $|a_\eps(x,t)|$. Due to inequality~\eqref{eq:ab} the \op\ $A^+_\eps$ is
 bounded, $A^+_\eps \le B_\eps$, and $\|A^+_\eps\| \le \|B_\eps\|$. Moreover,
 $A^+_\eps \cs A_\eps$ for any $\eps>0$ by construction. Since $B_\eps \to B$
 strongly as $\eps \to 0$, the norms $\|A^+_\eps\|$ are bounded uniformly
 in $\eps>0$ and the set $\{A^+_\eps\}_{\eps>0}$ is weakly compact in $\cB(\cH)$.
 Consequently there exists a sequence $\eps_n \to 0$ and an \op\ $A^+_0$ such
 that $A^+_{\eps_n} \to A^+_0$ weakly as $n \to \infty$. Passing to the limit
 $\eps_n \to 0$ in $A^+_\eps \le B_\eps$ and
 $A^+_\eps \cs A_\eps$, we find that $A^+_0 \le B$ and $A^+_0 \cs A$.
 Since $B\in\cM(A)$ was arbitrary, these inequalities prove that $A^+_0$ is
 the weak limit of $A^+_\eps$ as $\eps \to 0$ and that $A^+_0$ is the desired
 absolute value $|A|$ of~$A$.
\end{proof}

Let $\chi_a$ denote the multiplication \op\ by the characteristic function of
the interval~$[a,1]$ and
\[
   \cA := \{ A \in \cB(\cH) \mid A \chi_a = \chi_a A \chi_a, \quad
   \forall a \in [0,1]\}.
\]
$\cA$ is clearly a weakly closed subalgebra of $\cB(\cH)$. Put
\[
   \cA_+ := \cA \cap \cB_+(\cH)
\]
and note that $|A| \in \cA_+$ for any $A \in \cB_\cM(\cH) \cap \cA$.

\begin{remark}\label{rem:volt}
  Suppose that $A \in \cA_+$ and $A_\eps = R_\eps A R_\eps$, $\eps>0$.
  Then $A_\eps$ is an integral \op\ with a continuous kernel $a_\eps(x,t)$
  such that $a_\eps(x,t) = 0$ for $t>x$. Therefore $A_\eps$ is a Volterra
  \op\ (see, e.g., \cite[Sect.~68]{RN} and $A$ is a strong limit of a
  sequence of Volterra \op s.
\end{remark}

The above remark implies the following result.

\begin{lemma}\label{lem:com}
 Suppose that $A \in \cA_+$; then $[S,A]:= SA - AS$ belongs to~$\cA_+$.
\end{lemma}

\begin{proof}
 It suffices to notice that $[S,A_\eps]$ is an integral \op\ with the
 nonnegative kernel $(x-t)a_\eps(x,t)$ and hence belongs to $\cA_+$.
\end{proof}


\section{Similarity of the operators $S$ and $T$}\label{sec:simil}

In this section, we shall study the question of existence of an \op\
$K \in \cB(\cH)$ satisfying equation~\eqref{eq:simil}. Denoting $[S,K]:=SK-KS$,
we rewrite~\eqref{eq:simil} as the equation
\[
      [S,K] + VK + V = 0
\]
and apply a modification of Friedrichs successive approximation
method \cite[Ch.~II.6]{Fri} to solve the latter.

\begin{lemma}\label{lem:eq}
 Suppose that $V\in \cB_\cM(\cH)\cap\cA$ and
 $[S,W] \in \cM(V)$ for some $W \in \cA_+$. Then the equation
 \[
    [S,K] = V
 \]
 has a solution $K \in \cB_\cM(\cH)\cap \cA$ such that $K \cp W$.
\end{lemma}

\begin{proof}
 Since $[S,W]\cs V$, we have $R_\eps[S,W]R_\eps \cs R_\eps VR_\eps =:V_\eps$.
 By Lemma~\ref{lem:com} $R_\eps S \le S R_\eps$, whence
 \[
   R_\eps [S,W] R_\eps = R_\eps SW R_\eps - R_\eps WS R_\eps
    \le S W_\eps - W_\eps S = [S,W_\eps]
 \]
 and $[S,W_\eps] \cs V_\eps$, where $W_\eps =R_\eps W R_\eps$. It follows
 that the kernels $v_\eps$ and $w_\eps$ of the \op s $V_\eps$ and $W_\eps$
 satisfy the inequality
 \[
     \left| \frac{v_\eps(x,t)}{x-t}\right| \le w_\eps(x,t);
 \]
 in particular, $V_\eps \in\cA$.
 Denote by $K^{(\eps)}$ the integral \op\ induced by the kernel
 $v_\eps(x,t)/(x-t)$. By the above inequality $K^{(\eps)}$ is a bounded \op\
 belonging to the algebra~$\cA$,
 $K^{(\eps)} \cp W_\eps$, and $[S,K^{(\eps)}] = V_\eps$. Therefore the norms
 $\|K^{(\eps)}\|$ are bounded uniformly in~$\eps>0$ and there exist a bounded
 \op\ $K$ and a sequence $\eps_n\to0$ such that $K^{(\eps_n)}$ converge weakly
 to $K$ as $n\to\infty$.  Passing to the limit $\eps_n \to 0$ in the
 relations
 \[
        SK^{(\eps)} - K^{(\eps)} S = V_\eps, \qquad K^{(\eps)} \cp W_\eps
 \]
 we find that $K$ solves the equation $[S,K]=V$ and satisfies $K \cp W$.
 The proof is complete.
\end{proof}

\begin{corollary}\label{cor:eq}
  Suppose that $U,V \in \cB_\cM(\cH)\cap\cA$ and
  $[S,W] \in \cM(V)$ for some $W \in \cA_+$. Then the equation
 \[
    [S,K] = V U
 \]
 has a solution $K \in \cB_\cM(\cH)\cap \cA$ such that $K \cp W|U|$.
\end{corollary}

\begin{proof}
 Since $[S,W|U|] = [S,W]|U| + W[S,|U|] \ge [S,W]|U| \cs V U$ by
 Lemma~\ref{lem:com} and the assumptions of the corollary, Lemma~\ref{lem:eq}
 applies with $W|U|$ and $VU$ instead of $W$ and $V$, and the claim follows.
\end{proof}

\begin{theorem}\label{thm:eq}
 Suppose that $V \in \cA$ and that there exists an \op~$W\in\cA_+$ such that
 $\spr (W) <1$ and $[S,W]\cs V$. Then the equation
 \begin{equation}\label{eq:abstr}
        [S,K] + VK + V = 0
 \end{equation}
 has a solution $K \in \cB_\cM(\cH)\cap\cA$.
\end{theorem}

\begin{proof}
 We shall seek for $K$ of the form
 \[
        K = \sum_{k=1}^\infty K_n,
 \]
 where $K_n$ are found recursively from the relations
 \[
       [S,K_n] = - VK_{n-1}, \qquad n\in \bN,
 \]
 with $K_0:=I$. In virtue of Lemma~\ref{lem:eq} and Corollary~\ref{cor:eq}
 we find successively $K_n \in \cA$, $n=1,2,\dots,$ such that
 $K_n \cp W|K_{n-1}| \cp W^n$. Therefore $\|K_n\| \le \|W^n\|$, which shows
 that the series
 \(
   \sum_{n=1}^\infty K_n
 \)
 converges absolutely in the uniform operator topology and its sum $K$
 satisfies the equality
 \[
        [S,K] = \sum_{n=1}^\infty [S,K_n] = -\sum_{n=1}^\infty VK_{n-1}
        = -VK - V.
 \]
 The theorem is proved.
\end{proof}

\begin{corollary}\label{cor:spr}
 If under the assumptions of Theorem~\ref{thm:eq} $\spr (W) < 1/2$, then
 equation~\eqref{eq:abstr} has a solution $K\in\cB_\cM(\cH)\cap\cA$ with
 $\spr(K) < 1$.
\end{corollary}

\begin{proof}
 Observe that the solution $K$ constructed in the proof of
 Theorem~\ref{thm:eq} satisfies the inequality $K \cp W(I-W)^{-1}$.
 Therefore $\spr(K) \le \spr\bigl(W(I-W)^{-1}\bigr)$, and it suffices to note
 that $\spr \bigl(W(I-W)^{-1}\bigr) < 1$ if $\spr (W) < 1/2$.
\end{proof}


\section{Proof of the main results}\label{sec:appl}

\begin{proofof}{Proof of Theorem~\ref{thm:main}.}
It suffices to notice that the assumptions of Theorem~\ref{thm:eq} and
Corollary~\ref{cor:spr} are satisfied for the integral operator $W$
with the kernel
\[
   w (x,t) := \frac{|v(x,t)|}{x-t}.
\]
\end{proofof}

\begin{theorem}\label{thm:toepl}
 Suppose that $w$ is a positive kernel on $[0,1]\times[0,1]$ such that
 $w(x,t)=0$ if $t>x$ and $w(x,t) \le q(x-t)$ if $x\ge t$, where
 $q\in L_1(0,1)$. Then the induced integral \op~$W$ is a Volterra \op\
 in~$\cH$ and $\|W\| \le \| q\|_1$.
\end{theorem}

\begin{proof}
 We prove first that $W$ is bounded in~$\cH$. In fact,
 \begin{align*}
    \int_0^x w(x,t)\,dt &\le \int_0^x q(x-t)\,dt =
        \int_0^x q(t)\,dt \le \|q\|_1\\
    \int_t^1 w(x,t)\,dx &\le \int_t^1 q(x-t)\,dx =
        \int_0^{1-t} q(x)\,dx \le \|q\|_1,
 \end{align*}
 whence $\|W\| \le \|q\|_1$ by the Schur test~\cite[Theorem~5.2]{HS}.
 Next, put
\[
w^{(m)} =\ \biggr\{
        {\begin{aligned} w(x,t) \quad &\text{if }\ q(x-t) \le m,\\
                           0\quad &\text{if }\ q(x-t) > m.
        \end{aligned}}
\]
Then the induced integral \op~$W^{(m)}$ is a Volterra \op\
in~$\cH$ (see, e.g., \cite[Sect.~68]{RN}) and
\[
        \| W - W^{(m)}\| \le \int_{q>m} q(t)\,dt \to 0
\]
as $m\to \infty$ in view of $q\in L_1(0,1)$ and the above arguments.
Therefore $W$ is a Volterra operator as well and the lemma is proved.
\end{proof}

Corollary~\ref{cor:toepl} now easily follows from Theorems~\ref{thm:main} and
\ref{thm:toepl}.

\section{Some generalisations}\label{sec:gen}

In this section, we comment on some straightforward generalisations of the
main results. Observe first that the arguments of Sections~\ref{sec:nonneg}
and \ref{sec:simil} work for the Banach space $L_p(a,b)$ with $-\infty \le
a<b\le\infty$ arbitrary and any $p\in[1,\infty)$. Next,
$S$ can be replaced by the multiplication operator by any increasing
function $\phi(x)$. The analogue of Theorem~\ref{thm:main} reads as follows.

\begin{theorem*}{Theorem 1$'$} \em
Suppose that $\phi$ strictly increases on $(a,b)$ and that the kernel
\[
        w(x,t):=\Bigl|\frac{v(x,t)}{\phi(x)-\phi(t)}\Bigr|
\]
induces a bounded integral operator in $L_p(a,b)$ of spectral radius less
than $1/2$. Then the operators $S$ of multiplication by $\phi$ and $T:=S+V$,
where
\[
        (Vf)(x):= \int_a^x v(x,t) f(t)\,dt,
\]
are similar in $L_p(a,b)$. In particular, $S$ and $T$ are similar if 
$w(x,t) \le q(x-t)$ for some $q\in L_1(a,b)$.
\end{theorem*}

Observe also that most results of the paper hold in an arbitrary Banach
lattice~$\cX$ of functions over $(a,b)$ \cite{Sch} provided the identity
operator in $\cX$ is the strong limit of integral operators with continuous
kernels.

{\em Remark. } After this paper was finished, M.~M.~Malamud drew our attention 
to his note~\cite{Ma2}, where the statements of Theorem~\ref{thm:main} (under the 
assumption $\spr(W)=0$) and Corollary~\ref{cor:toepl} were announced without 
proof.



\begin{thebibliography}{Fre}

\bibitem{Fri} Friedrichs, K.~O.
	{\em Perturbation of Spectra in Hilbert Space}.
	Amer. Math. Soc., Providence, RI, 1965.

\bibitem{Fa} Faddeev, L.~D.
        On Friedrichs model in the perturbation theory of continuous spectrum//
        {\em Trudy MIAN}. 1964. V.~73. P.~292-313 (Russian).

\bibitem{Ko} Kokholm, N.~J.
        Spectral analysis of perturbed multiplication operators
        occurring in polymerisation chemistry//
        {\em Proc. Roy. Soc. Edinburgh Sect. A}. 1989. V.~113. No.1--2.
        P.~119--148.

\bibitem{NT} Naboko S.~N. and Tretter C.
        Lyapunov stability of a perturbed multiplication operator//
        {\em Contributions to operator theory in spaces with an indefinite 
	metric. The Heinz Langer anniversary volume on the occasion of his 
	60th birthday}. 
	Basel: Birkh\"auser, 1998. Oper. Theory Adv. Appl. V.~106.
	P.~309--326.

\bibitem{Fre} Freeman, J.~M.
	Volterra operators similar to $J:f(x) \mapsto \int_0^x f(y)\,dy$//
	{\em Trans. Amer. Math. Soc.} 1965. V.~116. P.~181--192.

\bibitem{Ma} Malamud, M.~M.
	Similarity of Volterra operators and related problems in the theory
	of differential equations of fractional orders//
	{\em Trudy Moscov. Mat. Obshch.} 1994. V.~55. P.~73--148 (Russian);
	transl. in {\em Trans. Moscow Math. Soc.} 1994. V.~55. P.~57-122.

\bibitem{KR} Krein M.~G., Rutman M.~A.
        Linear operators leaving invariant a cone in a Banach space//
        {\em Uspekhi Mat. Nauk}. 1948. V.~3. No.~1. P.~3--95 (Russian);
        transl. in {\em Functional Analysis and Measure Theory}, Transl.
        Amer.~Math.~Soc. 1962. V.~10. P.~5--105.

\bibitem{St} Stein E.~M.
        {\em Singular Integrals and Differentiability Properties of
        Functions}. Princeton Univ. Press, Princeton, N.~J., 1970.

\bibitem{HS} Halmos, P.~R., Sunder V.~S.
        {\em Bounded Integral Operators on $L^2$ Spaces}.
        Springer-Verlag, Berlin-Heidelberg-New York, 1978.

\bibitem{RN} Riesz F. and Sz.-Nagy B.
        {\em Le\c cons D'analyse Fonctionnelle}.
        Akad\'emiai Kiad\'o, Budapest, 1972.

\bibitem{Sch} Schaefer H.~H.
	{\em Banach Lattices and Positive Operators}.
	Springer-Verlag, Berlin, 1974.

\bibitem{Ma2} Malamud M.~M.//
	 {\em Uspekhi Mat. Nauk}. 1977. V.~32. No.~5. P.~167. (Russian)
\end{thebibliography}
\end{document}